\documentclass{article}

\usepackage{arxiv}

\usepackage{bm}
\usepackage{subcaption}
\usepackage{amsmath}
\usepackage{amsfonts}
\usepackage{listings}

\usepackage[svgnames]{xcolor}
\lstdefinelanguage{XML}
{
  basicstyle=\ttfamily\footnotesize,
  morestring=[b]",
  moredelim=[s][\bfseries\color{Maroon}]{<}{\ },
  moredelim=[s][\bfseries\color{Maroon}]{</}{>},
  moredelim=[l][\bfseries\color{Maroon}]{/>},
  moredelim=[l][\bfseries\color{Maroon}]{>},
  morecomment=[s]{<?}{?>},
  morecomment=[s]{<!--}{-->},
  commentstyle=\color{DarkOliveGreen},
  stringstyle=\color{blue},
  identifierstyle=\color{red}
}

\usepackage{color}

\usepackage{pgfplots}
\usepackage{pgfplotstable}
\pgfplotsset{
compat=newest, % version compability
tick label style={font=\footnotesize}, % size of the tick labels
}

\usepackage{overpic}

\usepackage{tikz}
\usetikzlibrary{mindmap,shadows}
\usetikzlibrary{arrows,intersections}
\usepackage{smartdiagram}

\usepackage{booktabs}
\usepackage[]{units}

\title{An approach using the null space to implement Dirichlet and constraint boundary conditions into FEM}

\author{ Stefan Schoder \\
	Group of Aeroacoustics and Vibroacoustics, IGTE\\
	TU Graz\\
	Inffeldgasse 18, 8010 Graz \\
	\texttt{stefan.schoder@tugraz.at} \\
}

\renewcommand{\shorttitle}{openCFS Software}

\usepackage{hyperref}
\hypersetup{
pdftitle={openCFS Software},
pdfsubject={preprint},
pdfauthor={Schoder, Roppert},
pdfkeywords={FEM Sotware, open Source},
}

\begin{document}
\maketitle

\begin{abstract}
A handy technique for the Finite Element Method (FEM) is presented that uses the null space for the implementation of Dirichlet and constraint boundary conditions. The focus of this method is to present an illustrative approach to modeling boundary constraints within FEM simulations for teaching. It presents a consistent way of including the boundary terms in the forcing and constructing the field solution after solving the algebraic system of equations.
\end{abstract}

% keywords can be removed

%\section*{Nomenclature}
%
% \begin{tabbing}
%   XXXXXXX \= \kill% this line sets tab stop
%   $c_0$ \> Speed of sound (mean value) \\
% \end{tabbing}

\section{Introduction}

Generally, essential boundary conditions are implemented by a penalty or an elimination method. This report describes an implementation-efficient way to incorporate Dirichlet boundary conditions into the final assembled system of equations,

\begin{equation}
\bm M \bm{ \ddot v}(t) + \bm D \bm{ \dot v}(t) + \bm K \bm v(t) = \bm F(t) \ .
\end{equation}

This type of equation arises in time-dependent finite element approaches. In the next section, we form a concept of incorporating the boundary conditions by a particular type of elimination.

\section{Implementation}
First, we conclude that any linear relation between the nodal degrees of freedom in $\bm v(t)$ can be expressed in the linear form
\begin{equation}
\bm b^{\rm T}\bm v(t) = v_{DB} \ .
\end{equation}

All finite element methods frequently use such types of Dirichlet boundary conditions. For example, the linear constraint $v_{\rm 1} = v_{\rm 2}$ is enforced by

\begin{equation}
\bm b^{\rm T} = [1 \ -1 \ 0 \ \cdots \ 0 ]  \ , 
\end{equation}

and the boundary condition $v_{\rm 1} = v_{DB \rm1}$ by

\begin{equation}
\bm b^{\rm T} = [1 \ 0 \ \cdots \ 0 ]  \ . 
\end{equation}

A matrix relation enforces multiple linear inhomogeneous constraints

\begin{equation}
\bm B \bm v(t) = \bm v_{DB} \ ,
\end{equation}

where each row of $\bm B$ represents one linear constraint. In general, we consider a system with a total number of constrained and unconstrained degrees of freedom of $ n $ and enforce $n_{enf}$ linear constraints on the system. Thus, the matrix $\bm B$ of size $n_{enf} \times n$ enforces the linear constraints.

Furthermore, we can introduce a variable transformation for the degrees of freedom $\phi_{\rm i} = v_{\rm i} - v_{DB \rm i} $. The transformation holds for all degrees of freedom since we set $v_{DB \rm i} = 0$ for all homogeneous constrained and unconstrained degrees of freedom. Thus, the multiple linear inhomogeneous constraints are simply enforced by the transformed matrix relation

\begin{equation}
\bm B \bm \phi(t) = \bm 0 \ ,
\end{equation}

and we found a homogeneous system of equations. 

We introduce a new set of $n-n_{enf}$ coordinates $\bm w(t)$ according to the relation 

\begin{equation}
\bm \phi(t) =\bm C \bm w(t) \ ,
\end{equation}
where $\bm C$ is a $n \times n-n_{enf}$ matrix of which the columns span the null space of $\bm B$. Hence $\bm C$ fulfills the relation

\begin{equation}
\bm B \bm C = \bm 0 \ .
\end{equation}

These coordinates $\bm w(t)$ fulfill the relation 

\begin{equation}
\bm B \bm \phi(t) = \bm B \bm C \bm w(t) = \bm 0 \ 
\end{equation}

per definition. This means that any displacement depicted by $\bm w(t)$ produces a displacement $\bm \phi(t)$, respectively $\bm v(t)$, that satisfies the constraint $\bm B \bm v(t) = \bm v_{DB}$! Another way to interpret this property is that the columns of $\bm C$ define a subspace of $\mathcal{R}^n$ in which the constraints are always satisfied. Thus, we choose a new set of variables such that the displacement is restricted to this subspace.

The next step is to pose the system equations in the new coordinates $\bm w(t)$. Here, we have to note that the basis functions and the test function are transformed to the new subspace since we are performing a Galerkin finite element method. For the unconstrained dynamic system of equations, we obtain

\begin{eqnarray}
 \bm C^{\rm T}  \bm M \bm C \bm{ \ddot w}(t) + \bm C^{\rm T} \bm D \bm C \bm{ \dot w}(t) +\bm C^{\rm T} \bm K \bm C \bm w(t) =\\
 = \bm C^{\rm T} \left( \bm F(t) +   \bm M \bm B^{\rm T} \bm{\ddot{v}}_{DB} +  \bm D \bm B^{\rm T} \bm{\dot{v}}_{DB} +  \bm K \bm B^{\rm T} \bm v_{DB } \right) \ .
\end{eqnarray}

where the constrained (projected) matrices $\bm C^{\rm T}  \bm M \bm C$, $\bm C^{\rm T}  \bm D \bm C$, and $\bm C^{\rm T}  \bm K \bm C$  are of size $(n-n_{enf}) \times (n-n_{enf})$. This leads to the projected solution $\bm w(t)$. Thus, the original system of equations is transformed to

\begin{equation}
 \bar{\bm M} \bm{ \ddot w}(t) + \bar{\bm D} \bm{ \dot w}(t) +\bar{\bm K} \bm w(t)  = \bar{ \bm F}(t)  \ .
\end{equation}
Therefore, we conserved the ability to treat the system of equations with the same temporal numerical procedures as the equations of motion. The constrained solution can be obtained by

\begin{equation}
\bm v(t) = \bm C \bm w(t) - \bm B^{\rm T} \bm v_{DB} \ .
\end{equation}

This method guarantees high efficiency since the matrices $\bm C $ and $\bm B $ since only some coefficients differ from zero. Thus, the sparsity pattern can be exploited. Furthermore, in the case of a nonlinear system, the matrices $\bm C $ and $\bm B $ do not change. However, the system matrices change; hence, we must recompute the matrices in every iteration.

\section{Reaction forces}
Additionally, due to the elegant way of implementing the boundaries, we can compute the reaction forces in the constrained node locations by this method. As described in \cite{ramers}, we split the equations into a part of Dirichlet constraints and the unconstrained part. We lack knowledge about the reaction forces in the part where we know the displacement field. We only consider a static problem for simplicity, but the extension towards dynamics is straightforward. The equation 

\begin{equation}
\bm F_{reac} = \bm B \bm K \bm C \bm w + \bm B \bm K \bm B^{\rm T} \bm v_{DB}  \ 
\end{equation}

follows from the method described by Rammerstorfer in \cite{ramers} at page 15. The first part of the right-hand side is the contribution of the unknown nodal values, and the second part is the contribution of the known Dirichlet boundaries.

\section{Conclusions}

The utilization of the null space to incorporate boundary conditions within the Finite Element Method (FEM) offers an interesting approach for a general and fully automated handling of boundaries without extensive code generation and testing as currently done for the elimination technique. It comes with an additional feature, making it possible to directly incorporate linear constraints on the solution variable. By seamlessly integrating Dirichlet and constraint boundary conditions, the approach has the potential to tackle complex engineering \cite{schoder2020hybrid,schoder2020computational,schoder2020aeroacoustic,schoder2022aeroacoustic,tieghi2023machine,schoder2022error,maurerlehner2022aeroacoustic,schoder2023acoustic,schoder2023dataset,wurzinger2024experimental,Kraxberger2023Validated} and facilitate the development of robust and reliable solutions in practical applications \cite{Schoder2022openCFS}. 

\bibliographystyle{abbrv}
\bibliography{references}  %%% Uncomment this line and comment out the ``thebibliography'' section below to use the external .bib file (using bibtex) .

%%% Uncomment this section and comment out the \bibliography{references} line above to use inline references.
% \begin{thebibliography}{1}

% 	\bibitem{kour2014real}
% 	George Kour and Raid Saabne.
% 	\newblock Real-time segmentation of on-line handwritten arabic script.
% 	\newblock In {\em Frontiers in Handwriting Recognition (ICFHR), 2014 14th
% 			International Conference on}, pages 417--422. IEEE, 2014.

% 	\bibitem{kour2014fast}
% 	George Kour and Raid Saabne.
% 	\newblock Fast classification of handwritten on-line arabic characters.
% 	\newblock In {\em Soft Computing and Pattern Recognition (SoCPaR), 2014 6th
% 			International Conference of}, pages 312--318. IEEE, 2014.

% 	\bibitem{hadash2018estimate}
% 	Guy Hadash, Einat Kermany, Boaz Carmeli, Ofer Lavi, George Kour, and Alon
% 	Jacovi.
% 	\newblock Estimate and replace: A novel approach to integrating deep neural
% 	networks with existing applications.
% 	\newblock {\em arXiv preprint arXiv:1804.09028}, 2018.

% \end{thebibliography}

\end{document}